\documentclass{amsart}
\usepackage{amsmath}
\usepackage{amssymb}
\usepackage{latexsym}
\usepackage{amsbsy}
\usepackage[all]{xy}
\usepackage{amscd}
\usepackage[dvips]{graphicx}
\usepackage{color}
\usepackage{amsmath,amssymb,latexsym,cancel,rotating}
\usepackage{graphicx,amssymb,mathrsfs,amsmath,color,fancyhdr,amsthm}
\usepackage[all]{xy}
\usepackage{leftidx}

\newtheorem{theorem}{Theorem}[section]
\newtheorem{lemma}[theorem]{Lemma}
\newtheorem{corollary}[theorem]{Corollary}
\newtheorem{definition}[theorem]{Definition}
\newtheorem{proposition}[theorem]{Proposition}

\newtheorem{example}[theorem]{Example}

\numberwithin{equation}{section}

\begin{document}

\title[On bases of quantum affine algebras]{On bases of quantum affine algebras}

\author[J. Xiao]{Jie Xiao}
\address{Department of Mathematics, Tsinghua University, Beijing 100084, P. R. China}
\email{jxiao@tsinghua.edu.cn (J. Xiao)}
\author[H. Xu]{Han Xu$^\ast$}
\address{Department of Mathematics, Tsinghua University, Beijing 100084, P. R. China}
\email{xu-h15@mails.tsinghua.edu.cn (H. Xu)}
\author[M. Zhao]{Minghui Zhao}
\address{School of Science, Beijing Forestry University, Beijing 100083, P. R. China}
\email{zhaomh@bjfu.edu.cn (M. Zhao)}

\thanks{This review is based on the report of Jie Xiao given in the conference "\emph{Forty Years of
Algebraic Groups, Algebraic Geometry, and Representation Theory in China}" in Jan. 5th, 2020}
\thanks{$^\ast$ supported by Tsinghua University Initiative Scientific Research Program (2019Z07L01006)}

%supported by Tsinghua University Initiative Scientific Research Program (2019Z07L01006)

%\institute{J. Xiao \at
%              Department of Mathematics, Tsinghua University, Beijing {\rm 100084}, P. R. China\\
%              \email{jxiao@tsinghua.edu.cn}
%           \and
%           H. Xu \at
%              Department of Mathematics, Tsinghua University, Beijing {\rm 100084}, P. R. China\\
%              \email{xu-h15@mails.tsinghua.edu.cn}
%            \and
%           M. Zhao \at
%              College of Science, Beijing Forestry University, Beijing {\rm 100083}, P. R. China\\
%              \email{zhaomh@bjfu.edu.cn}
%}

\date{\today}

%
%\keywords{}

%\begin{abstract}
%\end{abstract}

\bibliographystyle{abbrv}

\maketitle
\begin{center}
\emph{\text{In Memory of Prof. Xihua Cao}}
\end{center}

\setcounter{tocdepth}{1}
\tableofcontents

\section{Introduction}

Let $\textbf{U}^+$ be the positive part of the quantum group $\textbf{U}$ associated with a generalized Cartan matrix $A$.

In the case of finite type, Lusztig constructed the canonical basis $\textbf{B}$ of $\textbf{U}^+$ via two approaches (\cite{Lusztig_Canonical_bases_arising_from_quantized_enveloping_algebra}).
The first one is an elementary algebraic construction via Ringel-Hall algebra realization of $\textbf{U}^+$.
The isomorphism classes of representations of the corresponding Dynkin quiver form a PBW-type basis of $\textbf{U}^+$.
By a lemma (Lemma \ref{Lemma_Lusztig}) of Lusztig, one can construct a bar-invariant basis, which is the canonical basis $\textbf{B}$.
A remarkable characteristic in his construction is that Lusztig reveals the triangular relations among three kind of bases: PBW-basis, monomial basis and canonical basis.

The second one is a geometric construction.
In \cite{Lusztig_Canonical_bases_arising_from_quantized_enveloping_algebra} and \cite{Lusztig_Quivers_perverse_sheaves_and_the_quantized_enveloping_algebras}, Lusztig gave a geometric realization of $\mathbf{U}^+$ via the category of some semisimple complexes on the variety $E_{\nu}$ consisting of representations with dimension vector $\nu\in\mathbb{N}I$ of the corresponding quiver $Q$.
The set of the isomorphism classes of simple perverse sheaves gives the canonical basis $\mathbf{B}$ of $\mathbf{U}^+$.

The geometric construction of canonical basis was generalized to the cases of all types in \cite{Lusztig_Quivers_perverse_sheaves_and_the_quantized_enveloping_algebras} (see also \cite{Lusztig_Introduction_to_quantum_groups}). Furthermore, Lusztig in \cite{Lusztig_Affine_quivers_and_canonical_bases} gave the construction of affine canonical bases by the perverse sheaves associated with tame quivers, an important feature is that those perverse sheaves are indexed by the classes of aperiodic modules of tame quivers and irreducible modules of symmetric groups.

The generalization of this elementary algebraic construction to affine type is an important problem.

In \cite{BCP} and \cite{2002Beck-Nakajima-Crystal}, Beck-Chari-Pressley and Beck-Nakajima  defined a PBW-type basis and gave an algebraic construction of canonical basis. However, the meaning of this construction is not clear in terms of representation theory of quivers, in particular, it is not clear  how to parametrize the canonical basis in terms of aperiodic modules of tame quivers and irreducible modules of symmetric groups.

In the case of Kronecker quiver, Chen (\cite{Chen-Rootvector}) and Zhang (\cite{Zhang2000PBW}) defined a PBW-type basis by using Ringel-Hall algebra realization of $\textbf{U}^+$, McGerty (\cite{Mcgerty2004The_Kronecker_quiver}) gave the interpretation of the PBW-type basis of Beck-Nakajima in terms of representation theory of the Kronecker quiver.
In the case of cyclic quiver, Deng-Du-Xiao (\cite{DengDuXiao2007Generic}) defined a PBW-type basis and gave a concrete algebraic construction of the canonical basis by using the triangular relations between the monomial basis and the PBW-basis.

Lin-Xiao-Zhang in \cite{Lin_Xiao_Zhang_Representations_of_tame_quivers_and_affine_canonical_bases} provided a process to construct a PBW-type basis of $\textbf{U}^+$ and the canonical basis $\textbf{B}$ by using Ringel-Hall algebra approach. Recently
Xiao-Xu-Zhao  (\cite{Xiao-Xu-Zhao}) provided a direct method to construct a PBW-type basis of $\textbf{U}^+$ and the canonical basis $\textbf{B}$. Compared with the construction Lin-Xiao-Zhang in \cite{Lin_Xiao_Zhang_Representations_of_tame_quivers_and_affine_canonical_bases}, the PBW-type basis of Xiao-Xu-Zhao is a $\mathbb{Z}[v,v^{-1}]$-basis. Particularly the parametrization of the basis by aperiodic modules of tame quivers and irreducible characters of $S_n$ naturally arises in the construction.

In this paper, we shall review these constructions. Since we shall use representations of quivers, we mainly consider the quantum group corresponding to a symmetric generalized Cartan matrix.

\section{Preliminaries}

Let $I$ be a finite index set, $A=(a_{ij})_{i,j\in I}$ be a
symmetric
generalized Cartan matrix.
Denote by $\mathbf{U}$ the quantum group associated with the Cartan matrix $A$ and $\mathbf{U}^+$ the positive part of $\mathbf{U}$ (\cite{Lusztig_Introduction_to_quantum_groups}).

Fix an indeterminate $v$. For any $n\in \mathbb{Z}$, set
$$[n]_{v}=\frac{v^n-v^{-n}}{v-v^{-1}}.$$
Let $[0]_{v}!=1$ and $[n]_{v}!=[n]_{v}[n-1]_{v}\cdots[1]_{v}$ for any $n\in\mathbb{Z}_{>0}$.

Note that $\mathbf{U}^+$ is an associate algebra over $\mathbb{Q}(v)$ generated by
the elements $E_i$ for various $i\in I$ subject to the quantum Serre relations
$$\sum_{k=0}^{1-a_{ij}}(-1)^{k}E_i^{(k)}E_jE_i^{(1-a_{ij}-k)}=0$$ for all $i\neq j\in I$, where $E_i^{(n)}=E_i^n/[n]_{v}!$.

Let $\,\,\bar{}:\mathbf{U}^+\rightarrow\mathbf{U}^+$ be the unique $\mathbb{Q}$-algebra involution such that $$\overline{v^{n}}=v^{-n}\,\,\,\,\,\,\textrm{ and }\,\,\,\,\,\,\overline{E_i}=E_i.$$

Let $Q=(I,H,s,t)$ be a quiver  without loops, where $I$ is the set of vertices, $H$ is the set of arrows and $s,t:H\rightarrow I$ are two maps sending an arrow $h\in H$ to the source $s(h)$ and  target $t(h)$ respectively.
Let $$a_{ij}=|\{h\in H\,\,|\,\,\textrm{$s(h)=i$ and $t(h)=j$}\}|+|\{h\in H\,\,|\,\,\textrm{$s(h)=j$ and $t(h)=i$}\}|.$$
Then $A=(a_{ij})_{i,j\in I}$ is a symmetric generalized Cartan matrix and called the generalized Cartan matrix associated to the quiver $Q$.
The Euler form on $\mathbb{Z}I$ is defined as
$$\langle\nu,\nu'\rangle=\sum_{i\in I}\nu_i\nu'_i-\sum_{h\in H}\nu_{s(h)}\nu'_{t(h)}$$
and the symmetric Euler form is defined as
$(\nu,\nu')=\langle\nu,\nu'\rangle+\langle\nu',\nu\rangle$
for any  $\nu=\sum_{i\in I}\nu_ii$ and $\nu'=\sum_{i\in I}\nu'_ii$.

Let $k=\mathbb{F}_q$ be a finite field with $q$ elements and $kQ$ be the path algebra.
Denote by $\textrm{rep}_{k}Q$ the abelian category of finite dimensional representations of $Q$ over $k$.
There is an  isomorphism between $\textrm{rep}_{k}Q$ and the category  $\textrm{mod-}kQ$ of finite dimensional  $kQ$-modules.

Let $\mathcal{P}_k$ be the set of isomorphism classes of objects in $\textrm{rep}_{k}Q$.
For any $\alpha\in\mathcal{P}_k$, choose an object $M_{\alpha}\in\textrm{rep}_{k}Q$ such that the isomorphism classes $[M_{\alpha}]=\alpha$.
Define the dimension vector of $\alpha\in\mathcal{P}_k$ by $\underline{\dim}\,\alpha=\underline{\dim} M_{\alpha}$.

Given three elements $\alpha_1,\alpha_2$ and $\alpha$ in $\mathcal{P}_k$, denote by $g^{\alpha}_{\alpha_1\alpha_2}$ the number of subrepresentations $N$ of $M_{\alpha}$ such that $N\simeq M_{\alpha_2}$ and $M_{\alpha}/N\simeq M_{\alpha_1}$ in $\textrm{rep}_{k}Q$.
Let $v_q=\sqrt{q}\in \mathbb{C}$.  The twisted Ringel-Hall algebra $\mathcal{H}^{\ast}(kQ)$ is the $\mathbb{Q}(v_q)$-space with basis $$\{u_{\alpha}\,\,|\,\,\alpha\in\mathcal{P}_k\},$$whose multiplication is given by
\begin{displaymath}
u_{\alpha_1}\ast u_{\alpha_2}=\sum_{\alpha\in\mathcal{P}} v_q^{\langle\underline{\dim}\alpha_1,\underline{\dim}\alpha_2\rangle}g^{\alpha}_{\alpha_1,\alpha_2}u_{\alpha}.
\end{displaymath}

For any $\alpha\in\mathcal{P}_k$, let $$\langle{M_{\alpha}}\rangle=v_q^{-\dim M_{\alpha}+\dim\textrm{End}M_{\alpha}}u_{\alpha}.$$
The set $\{\langle{M_{\alpha}}\rangle\,\,|\,\,\alpha\in\mathcal{P}_k\}$ is also a $\mathbb{Q}(v_q)$-basis of $\mathcal{H}^{\ast}(kQ)$.
There is a bilinear form $(-,-)$ on $\mathcal{H}^{\ast}(kQ)$ defined in \cite{green1995hall}.

Denote by  $\mathcal{C}^{\ast}(kQ)$ the composition subalgebra of $\mathcal{H}^{\ast}(kQ)$ generated by $u_{i}=u_{[S_i]}$ for all $i\in I$, where $S_i$ is the simple $kQ$-module
corresponding to $i\in I$.

Then we shall recall the generic form of $\mathcal{C}^{\ast}(kQ)$.
Let $\mathcal{K}$ be a set of some finite fields $k$ such that the set $\{q_k=|k|\,\,|\,\,k\in\mathcal{K}\}$ is an infinite set.
Consider the direct product
\begin{displaymath}
\mathcal{H}^{\ast}(Q)=\prod_{k\in\mathcal{K}}\mathcal{H}^{\ast}(kQ)
\end{displaymath}
and the elements $v=(v_{q_k})_{k\in\mathcal{K}}$, $v^{-1}=(v_{q_k}^{-1})_{k\in\mathcal{K}}$ and $u_i=(u_i(k))_{k\in\mathcal{K}}$.
By $\mathcal{C}^{\ast}(Q)_{\mathbb{Q}[v,v^{-1}]}$ we denote the subalgebra of $\mathcal{H}^{\ast}(Q)$ generated by $v$, $v^{-1}$ and $u_i$ over $\mathbb{Q}$. We may regard it as a $\mathbb{Q}[v,v^{-1}]$-algebra generated by $u_i$, where $v$ is viewed as an indeterminate. Finally, define $\mathcal{C}^{\ast}(Q)=\mathbb{Q}(v)\otimes_{\mathbb{Q}[v,v^{-1}]}\mathcal{C}^{\ast}(Q)_{\mathbb{Q}[v,v^{-1}]}$, which is called the generic twisted composition algebra of $Q$.

\begin{theorem}[\cite{green1995hall}\cite{Ringel_Hall_algebras_and_quantum_groups}]
Let $Q$ be a quiver without loops, $A$ the corresponding generalized Cartan matrix and $\mathbf{U}^+$ the positive part of quantum group of type $A$.
There is an isomorphism of $\mathbb{Q}(v)$-algebras:
\begin{eqnarray*}
\mathcal{C}^{\ast}(Q)&\cong&\mathbf{U}^+\\
u_i&\mapsto&E_i.
\end{eqnarray*}
\end{theorem}

Under this isomorphism, the bar involution on $\mathbf{U}^+$ induces a bar involution
$\bar{}:\mathcal{C}^{\ast}(Q)\rightarrow\mathcal{C}^{\ast}(Q)$ such that $\overline{v^{n}}=v^{-n}$ and $\overline{u_i}=u_i.$

In \cite{Lusztig_Canonical_bases_arising_from_quantized_enveloping_algebra,Lusztig_Quivers_perverse_sheaves_and_the_quantized_enveloping_algebras}, Lusztig gave a geometric realization of $\mathbf{U}^+$. Let $Q$ be a quiver without
loops with associated generalized Cartan matrix $A$.
Lusztig considered the variety $E_{\nu}$ consisting of representations with dimension vector $\nu\in\mathbb{N}I$ of the quiver $Q$ over algebraically
closed field $\bar{k}$ and the category $\mathcal{Q}_{\nu}$ of some semisimple complexes of constructible sheaves on $E_{\nu}$.
Let $K(\mathcal{Q}_{\nu})$ be the Grothendieck group of $\mathcal{Q}_{\nu}$.
Considering all dimension vectors, let
$$K(\mathcal{Q})=\bigoplus_{\nu}K(\mathcal{Q}_{\nu}).$$
Lusztig define induction functors on $\mathcal{Q}_{\nu}$ and get a $\mathbb{Z}[v,v^{-1}]$-algebra structure on $K(\mathcal{Q})$.
This algebra is isomorphic to the integral form of $\mathbf{U}^+$.
The set $\mathbf{B}$ of the isomorphism classes of simple perverse sheaves gives a basis of $\mathbf{U}^+$, which is called the canonical basis.

\section{Canonical bases of finite types}

Assume that $Q$ is a Dynkin quiver. In this case, the algebraic construction of canonical basis was introduced by Lusztig
in \cite{Lusztig_Canonical_bases_arising_from_quantized_enveloping_algebra} (see also \cite{Rongel_The_Hall_algebra_approach_to_quantum_groups}).

Let $A=(a_{ij})_{i,j\in I}$ be the generalized Cartan matrix associated to the quiver $Q$ and
denote by $\Delta^+$ the set of positive roots of the Lie algebra $\mathfrak{g}(A)$ with $i$ corresponding to
simple roots. The dimension vector $\underline{\dim}$ induces a bijection between the set of isomorphism classes of indecomposable objects ind-$\mathcal{P}$ and the set $\Delta^+$ by Gabriel theorem. Given a positive root $\alpha$, choose an indecomposable representation $M_\alpha$ of $Q$ such that $[M_\alpha]=\alpha$.

Denote by $\mathbb{N}^{\Delta^+}$ the set of all functions $\phi:\Delta^+\rightarrow\mathbb{N}$. For each $\phi:\Delta^+\rightarrow\mathbb{N}$, define a representation
\begin{displaymath}
M_{\phi}=\bigoplus_{\alpha\in\Delta^+}M_\alpha^{\oplus\phi(\alpha)}.
\end{displaymath}
Then $\mathcal{P}=\{[M_{\phi}]\,\,|\,\,\phi\in\mathbb{N}^{\Delta^+}\}$.

Now $\mathcal{H}^{\ast}(kQ)$ is spanned by the set $$\{u_{\phi}=u_{[M_{\phi}]}\,\,|\,\,\phi:\Delta^+\rightarrow\mathbb{N}\}$$ as $\mathbb{Q}(v)$-vector space, which is called a PBW-type basis.

Since $Q$ is representation-directed, we can define a total order on the set $\Delta^+$
such that $$\textrm{Hom}(M_{\alpha},M_{\beta})\not=0\Rightarrow\alpha\leq\beta$$ for any $\alpha,\beta\in\Delta^+$.
This total order induces an order on $\mathbb{N}^{\Delta^+}$.
For any $\phi,\psi:\Delta^+\rightarrow\mathbb{N}$, define $\phi<\psi$ if and only if there exists $\alpha\in\Delta^+$ such that $\phi(\alpha)>\psi(\alpha)$ and $\phi(\beta)=\psi(\beta)$ for all $\alpha>\beta\in\Delta^+$.

For each $\phi:\Delta^+\rightarrow\mathbb{N}$, there exists a monomial $m_{\phi}$ on the divided powers of Chevalley generators $u_{i}$  satisfying
\begin{displaymath}
m_{\phi}=\langle{M_{\phi}}\rangle+\sum_{\phi'<\phi}a^{\phi'}_{\phi}\langle{M_{\phi'}}\rangle,
\end{displaymath}
with $a^{\phi'}_{\phi}\in \mathbb{Z}[v,v^{-1}]$. Since $\overline{m_{\phi}}=m_{\phi}$, we have
\begin{displaymath}
\overline{\langle{M_{\phi}}\rangle}=\langle{M_{\phi}}\rangle+\sum_{\phi'<\phi}b^{\phi'}_{\phi}\langle{M_{\phi'}}\rangle,
\end{displaymath}
with $b^{\phi'}_{\phi}\in \mathbb{Z}[v,v^{-1}]$ such that
\begin{enumerate}
  \item[(1)]$b^{\phi}_{\phi}=1$ for all $\phi$ in $\mathbb{N}^{\Delta^+}$;
  \item[(2)]for all $\phi'\leq \phi$ in $\mathbb{N}^{\Delta^+}$,$$\sum_{\phi'',\phi'\leq \phi''\leq \phi}\overline{b^{\phi'}_{\phi''}}b^{\phi''}_{\phi}=\delta_{\phi,\phi'}.$$
\end{enumerate}

Here we need a lemma by Lusztig, which can be obtained by an elementary linear algebra method.

\begin{lemma}[\cite{Lusztig_Introduction_to_quantum_groups}]\label{Lemma_Lusztig}
Let $H$ be a set with a partial order $\leq$ such that for any $h'\leq h$ in $H$, the set $\{h''\,\,|\,\,h'\leq h''\leq h\}$ is finite. Assume that for each $h'\leq h$ in $H$, we are given an element $r_{h}^{h'}\in\mathbb{Z}[v,v^{-1}]$ such that
\begin{enumerate}
  \item[(1)]$r_{h}^{h}=1$ for all $h$ in $H$;
  \item[(2)]for all $h'\leq h$ in $H$,$$\sum_{h'',h'\leq h''\leq h}\overline{r^{h'}_{h''}}r^{h''}_{h}=\delta_{h,h'}.$$
\end{enumerate}
Then there is a unique family of elements $p_{h}^{h'}\in\mathbb{Z}[v^{-1}]$ defined for all $h'\leq h$ in $H$ such that
\begin{enumerate}
  \item[(1)]$p_{h}^{h}=1$ for all $h$ in $H$;
  \item[(2)]$p_{h}^{h'}\in v^{-1}\mathbb{Z}[v^{-1}]$ for all $h'\leq h$ in $H$;
  \item[(3)]for all $h'\leq h$ in $H$, $$p_{h}^{h'}=\sum_{h'',h'\leq h''\leq h}\overline{p_{h''}^{h'}}r_{h}^{h''}.$$
\end{enumerate}
\end{lemma}

By Lemma \ref{Lemma_Lusztig}, there exists a unique family of elements $c^{\phi'}_{\phi}\in\mathbb{Z}[v^{-1}]$ defined for all $\phi'\leq\phi$ in $\mathbb{N}^{\Delta^+}$ such that
\begin{enumerate}
  \item[(1)]$c^{\phi}_{\phi}=1$ for all $\phi$ in $\mathbb{N}^{\Delta^+}$;
  \item[(2)]$c^{\phi'}_{\phi}\in v^{-1}\mathbb{Z}[v^{-1}]$ for all $\phi'\leq\phi$ in $\mathbb{N}^{\Delta^+}$;
  \item[(3)]for all $\phi'\leq\phi$ in $\mathbb{N}^{\Delta^+}$, $$c^{\phi'}_{\phi}=\sum_{\phi'',\phi'\leq \phi''\leq \phi}\overline{c^{\phi'}_{\phi''}}b^{\phi''}_{\phi}.$$
\end{enumerate}

For any $\phi\in\mathbb{N}^{\Delta^+}$, let
$$C_{\phi}=\langle{M_{\phi}}\rangle+\sum_{\phi'<\phi}c^{\phi'}_{\phi}\langle{M_{\phi'}}\rangle.$$
These formulas hold for every finite fields and may be viewed as formulas in $\mathcal{H}^{\ast}(Q)=\mathcal{C}^{\ast}(Q)$.
Then we have the following theorem.
\begin{theorem}\label{canonical-finite}
The set $\{C_{\phi}\,\,|\,\,\phi:\Delta^+\rightarrow\mathbb{N}\}$ is a $\mathbb{Z}[v,v^{-1}]$-basis of $\mathcal{H}^{\ast}(Q)_{\mathbb{Z}[v,v^{-1}]}$ satisfying the following conditions.
\begin{enumerate}
  \item[(1)]$\overline{C_{\phi}}=C_{\phi}$;
  \item[(2)]$(C_{\phi},C_{\phi'})\in\delta_{\phi,\phi'}+v^{-1}\mathbb{Z}[[v^{-1}]]\cap\mathbb{Q}(v)$.
\end{enumerate}
\end{theorem}

Under the isomorphism between $\mathcal{H}^{\ast}(Q)$ and $\mathbf{U}^+$, the set $$\{C_{\phi}\,\,|\,\,\phi:\Delta^+\rightarrow\mathbb{N}\}$$ induces a basis of $\mathbf{U}^+$. This basis is just the canonical basis $\mathbf{B}$ of $\mathbf{U}^+$, by Theorem \ref{canonical-finite} and the uniqueness of canonical basis of $\mathbf{U}^+$.

\begin{example}
Take the quiver $Q$ of type $A_3$ for example.
$$Q: \xymatrix{
1\ar[r]&2\ar[r]&3
}.$$
The AR-quiver is as following.
$$
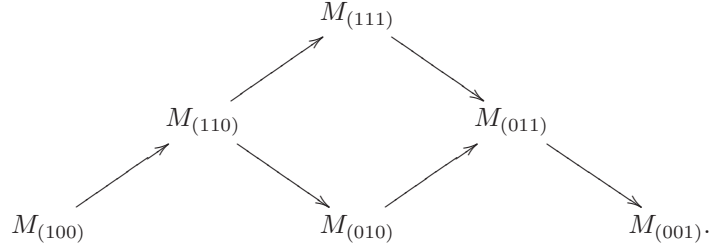
\xymatrix{
&&M_{(111)}\ar[rd]&&\\
&M_{(110)}\ar[ur]\ar[dr]&&M_{(011)}\ar[rd]&\\
M_{(100)}\ar[ur]&&M_{(010)}\ar[ru]& &M_{(001)}
.}$$

For dimension vector $\nu=(111)$, there exist isomorphism classes of the following modules
$$[M_{(111)}],[M_{(110)}\oplus M_{(001)}],[M_{(100)}\oplus M_{(011)}],[M_{(100)}\oplus M_{(010)}\oplus M_{(001)}].$$
Hence,
$$\{\langle{M_{(111)}}\rangle,\langle{M_{(110)}\oplus M_{(001)}}\rangle,\langle{M_{(100)}\oplus M_{(011)}}\rangle,\langle{M_{(100)}\oplus M_{(010)}\oplus M_{(001)}}\rangle\}$$
are the elements in the PBW-type basis with dimension vector $\nu=(111)$.

By construction, the elements in the corresponding monomial basis are
\begin{eqnarray*}
u_1u_2u_3=&\langle{M_{(111)}}\rangle+v^{-2}\langle{M_{(110)}\oplus M_{(001)}}\rangle+v^{-2}\langle{M_{(100)}\oplus M_{(011)}}\rangle\\
&+v^{-3}\langle{M_{(100)}\oplus M_{(010)}\oplus M_{(001)}}\rangle,
\end{eqnarray*}
$$u_3u_1u_2=\langle{M_{(110)}\oplus M_{(001)}}\rangle+v^{-1}\langle{M_{(100)}\oplus M_{(010)}\oplus M_{(001)}}\rangle,$$
$$u_2u_3u_1=\langle{M_{(100)}\oplus M_{(011)}}\rangle+v^{-1}\langle{M_{(100)}\oplus M_{(010)}\oplus M_{(001)}}\rangle,$$
$$u_3u_2u_1=\langle{M_{(100)}\oplus M_{(010)}\oplus M_{(001)}}\rangle,$$
which are also the elements with dimension vector $(111)$ in the canonical basis.

It is clear that these PBW-type basis elements are the leading terms of the corresponding canonical basis elements.

\end{example}

\section{Beck-Nakajima's construction}\label{BN}

In this section, we shall recall the construction of canonical basis given by Beck-Nakajima (\cite{2002Beck-Nakajima-Crystal}\cite{Nakajima_Crystal_canonical_and_PBW_bases_of_quantum_affine_algebras}).

Let $A=(a_{ij})_{i,j\in I}$ be a generalized Cartan matrix of affine type, where $I=\{0,1,\ldots,n\}$, $0\in I$ is the exceptional point and $I_0=I\backslash\{0\}$. Let $$D=\textrm{diag}(d_0,d_1,\ldots,d_n)$$ be a diagonal matrix such that $DA$ is symmetric.
Let $\Delta^{+}$ be the set of positive roots and $R$ the set of all positive real roots. Let $\{\alpha_i\,\,|\,\,i\in I\}$ be the set of simple roots. Let $v_i=v^{d_i}$.
%{\bc Note that the positive part $\mathbf{U}^+$ is an associate algebra over $\mathbb{Q}(v)$ generated by
%the elements $E_i$ for various $i\in I$ subject to the quantum Serre relations
%$$\sum_{k=0}^{1-a_{ij}}(-1)^{k}E_i^{(k)}E_jE_i^{(1-a_{ij}-k)}=0$$ for all $i\neq j\in I$, where $E_i^{(n)}=E_i^n/[n]_{v_i}!$.}

We follow the notations of \cite{2002Beck-Nakajima-Crystal}. Denote by $\hat{W}$ the affine Weyl group generated by simple reflections $s_i$ for $i\in I$.
% where $$ s_i(\lambda)=\lambda-(\lambda,\hat{\alpha}_i)\alpha_i.$$
Let $\tilde{W}$ be the extended affine Weyl group. Then $\tilde{W}=J\ltimes\hat{W}$, where $J$
%$J=\{w\in\tilde{W}\,\,|\,\,w(\Delta^{+})\subset\Delta^{+}\}$
is a subgroup of the group of Dynkin diagram automorphism and $\tau s_i=s_{\tau(i)}\tau\in\tilde{W}$ for $\tau\in J, s_i\in\hat{W}$. For any $i\in I_0$, there exists $t_{\tilde{w}_i}\in\tilde{W}$ such that $$
t_{\tilde{w}_i}(\alpha_j)=\left\{
\begin{array}{cc}
   \alpha_j & \textrm{if $j\neq i,0$}\\
   \alpha_i-d_i\delta & \textrm{if $j=i$}\\
   \alpha_0+a_id_i\delta & \textrm{if $j=0$},
\end{array}
\right.$$
where the minimal imaginary positive root $\delta=\sum a_i\alpha_i$ and $a_0=1$.

Let $s_{i_1}s_{i_2}\cdots s_{i_N}\tau$ be a reduced expression of $t_{\tilde{w}_n}t_{\tilde{w}_{n-1}}\cdots t_{\tilde{w}_1}$. Define an infinite sequence $$h=(\cdots,i_{-1},i_0,i_1,\cdots)$$
in $I$ such that $i_{k+N}=\tau(i_k)$ for any $k\in\mathbb{Z}$.
Let $$R_{<}=\{\beta_0=\alpha_{i_0},\beta_{-1}=s_{i_0}(\alpha_{i_{-1}}),\beta_{-2}=s_{i_0}s_{i_{-1}}(\alpha_{i_{-2}}),\cdots\}$$
and
$$R_{>}=\{\beta_{1}=\alpha_{i_1},\beta_{2}=s_{i_1}(\alpha_{i_2}),\beta_{3}=s_{i_1}s_{i_2}(\alpha_{i_2}),\cdots\}.$$
It is well-known that $$R=R_{>}\bigsqcup R_{<}.$$

For all $j\in I$, denote by $T_j$ the Lusztig's symmetries $T''_{j,1}$ in \cite{Lusztig_Introduction_to_quantum_groups}. For any $k\in\mathbb{Z}_{>0}$, let $$E_{\beta_k}=T_{i_1}T_{i_2}\cdots T_{i_{k-1}}(E_{i_k}).$$
For any $k\in\mathbb{Z}_{\leq0}$, let $$E_{\beta_k}=T^{-1}_{i_0}T^{-1}_{i_{-1}}\cdots T^{-1}_{i_{k+1}}(E_{i_k}).$$
Then $E_{\beta_k}$ are the root vectors for the real roots $\beta_k\in R$.

Then we shall define imaginary root vectors. For $k>0$ and $i\in I_0$, let
$$\tilde{\Psi}_{i,kd_i}=E_{kd_i\delta-\alpha_i}E_{\alpha_i}-v^{-2}_{i}E_{\alpha_i}E_{kd_i\delta-\alpha_i},$$
$\tilde{P}_{i,0}=1$ and
$$\tilde{P}_{i,kd_i}=\left\{\begin{array}{cc}
                             \frac{1}{[2k]_n}\sum_{s=1}^{k}v^{2(s-k)}_{n}\tilde{\Psi}_{n,s}\tilde{P}_{n,k-s} & \textrm{if $i=n$ and $A$ is of type $A^{(2)}_{2n}$}\\
                             \frac{1}{[k]_i}\sum_{s=1}^{k}v^{s-k}_{i}\tilde{\Psi}_{i,sd_i}\tilde{P}_{i,(k-s)d_i} & \textrm{otherwise}.\\
                           \end{array}
\right.$$

Let $\mathscr{P}$ be the set of all partitions and $c_0:I_0\rightarrow \mathscr{P}$ be a map.
For $\lambda=(\lambda_1\geq \lambda_2\geq\cdots)\in\mathscr{P}$, define $$S_{\lambda}=\det(\tilde{P}_{i,(\lambda_k-k+m)d_i})_{1\leq k,m\leq t},$$
where $t$ is the length of $\lambda$.
% and ${p^{(i)}}^t$ is the transpose of $p^{(i)}$.
Denote $$S_{c_0}=\prod_{i=1}^{n}S_{c_0(i)}.$$

Let $\mathcal{E}$ be the set of all such $\bar{c}=(c,c_0)$, where $c_0:I_0\rightarrow \mathscr{P}$ is a map and $c:\mathbb{Z}\rightarrow\mathbb{N}$ is a function with finite support.
For any $\bar{c}\in\mathcal{E}$ and $p\in\mathbb{Z}$, let
\begin{eqnarray*}
L(\bar{c},p)&=&\left(E_{i_p}^{(c(p))}T_{i_p}^{-1}(E_{i_{p-1}}^{(c(p-1))})T_{i_p}^{-1}T_{i_{p-1}}^{-1}(E_{i_{p-2}}^{(c(p-2))})\cdots\right)\\
&\times&T_{i_{p+1}}T_{i_{p+2}}\cdots T_{i_0}(S_{c_0})\\
&\times&\left(\cdots T_{i_{p+1}}T_{i_{p+2}}(E_{i_{p+3}}^{(c(p+3))})T_{i_{p+1}}(E_{i_{p+2}}^{(c(p+2))})E_{i_{p+1}}^{(c(p+1))}\right),
\end{eqnarray*}
where $i_p$ are from the sequence $h$.

For any $p\in\mathbb{Z}$, Beck-Nakajima defined a partial ordering $<_{p}$ on $\mathcal{E}$ such that the following Theorem holds.

\begin{theorem}[\cite{2002Beck-Nakajima-Crystal}]The set $\{L(\bar{c},p)\,\,|\,\,\bar{c}\in\mathcal{E}, p\in\mathbb{Z}\}$ is a $\mathbb{Z}[v,v^{-1}]$-basis of $\mathbf{U}_{\mathbb{Z}[v,v^{-1}]}^{+}$ such that
\begin{enumerate}
  \item[(1)]$(L(\bar{c},p),L(\bar{c}',p))\in\delta_{\bar{c},\bar{c}'}+v^{-1}\mathbb{Z}[[v^{-1}]]\cap\mathbb{Q}(v)$;
  \item[(2)]$$\overline{L(\bar{c},p)}=L(\bar{c},p)+\sum_{\bar{c}'<_{p}\bar{c}}a_{\bar{c}\bar{c}'}L(\bar{c}',p)$$ with $a_{\bar{c}\bar{c}'}\in\mathbb{Q}(v)$.
\end{enumerate}
\end{theorem}

The set $\{L(\bar{c},p)\,\,|\,\,\bar{c}\in\mathcal{E}, p\in\mathbb{Z}\}$ is called a PBW-type basis of $\mathbf{U}^+$.

Beck-Nakajima also proved the following Theorem.
\begin{theorem}[\cite{2002Beck-Nakajima-Crystal}]
For any $\bar{c}\in\mathcal{E}$ and $p\in\mathbb{Z}$, there exists a unique $b(\bar{c},p)\in \mathbf{U}_{\mathbb{Z}[v,v^{-1}]}^{+}$ satisfying the following conditions
\begin{enumerate}
  \item[(1)]$\overline{b(\bar{c},p)}=b(\bar{c},p)$;
  \item[(2)]$(b(\bar{c},p),b(\bar{c}',p))\in\delta_{\bar{c},\bar{c}'}+v^{-1}\mathbb{Z}[[v^{-1}]]\cap\mathbb{Q}(v)$;
  \item[(3)]$$b(\bar{c},p)=L(\bar{c},p)+\sum_{\bar{c}'<_{p}\bar{c}}\xi_{\bar{c}\bar{c}'}L(\bar{c}',p)$$ with $\xi_{\bar{c}\bar{c}'}\in v^{-1}\mathbb{Z}[v^{-1}]$.
\end{enumerate}
Moreover, the set $\{b(\bar{c},p)\,\,|\,\,\bar{c}\in\mathcal{E}, p\in\mathbb{Z}\}$ is the canonical basis of $\mathbf{U}^{+}$.
\end{theorem}

\section{Kronecker quiver}\label{Kronecker}

Let $Q$ be the Kronecker quiver with $I=\{0,1\}$ and $H=\{\rho_1,\rho_2\}$:
$$\xymatrix@=8em{
0\ar@/^/[r]^{\rho_1}\ar@/_/[r]_{\rho_2} & 1
}$$
Let $kQ$ be the path algebra of $Q$ over finite field $k$.

The set of dimension vectors of indecomposable $kQ$-modules is
\begin{displaymath}
\Delta^{+}=\{(l+1,l),(m,m),(n,n+1)\,\,|\,\,l,m,n\in\mathbb{Z},l\geq 0,m\geq 1,n\geq 0\}.
\end{displaymath}
The dimension vectors $(l+1,l)$ and $(n,n+1)$ correspond to preprojective and preinjective indecomposable $kQ$-modules respectively.
The subset consisting of $(l+1,l)$ (resp. $(n,n+1)$) is denoted by $Prep$ (resp. $Prei$).

For any $n\in\mathbb{N}$, let $M(n+1,n)$ and $M(n,n+1)$ be the indecomposable $kQ$-module of dimension vectors $(n+1,n)$ and $(n,n+1)$ respectively.
For real root vectors, define
$$
E_{(n+1,n)}=\langle{M(n+1,n)}\rangle$$and$$E_{(n,n+1)}=\langle{M(n,n+1)}\rangle.
$$

Then we shall define imaginary root vectors. Let $\delta=(1,1)$. As a special case of the definition in Section \ref{BN}, let
$$\tilde{\Psi}_{k}=\tilde{\Psi}_{1,k}=E_{(k-1,k)}E_{(1,0)}-v^{-2}E_{(1,0)}E_{(k-1,k)},$$
$\tilde{P}_{0}=1$ and
$$\tilde{P}_{k}=\frac{1}{[k]}\sum_{s=1}^{k}v^{s-k}\tilde{\Psi}_{s}\tilde{P}_{k-s},$$
for $k>0$.

\begin{proposition}[\cite{Chen-Rootvector}\cite{Zhang2000PBW}\cite{Mcgerty2004The_Kronecker_quiver}]
It holds that
$$\tilde{P}_{k}=\sum_{[M]:\underline{\dim}M=k\delta\atop\textrm{$M$ is regular}}v^{-\dim M}u_{[M]}$$
for $k\in\mathbb{Z}_{>0}$.
\end{proposition}

For any partition $\lambda=(\lambda_1\geq\lambda_2\geq\cdots\geq\lambda_t)$, let
\begin{displaymath}
\tilde{P}_{\lambda}=\tilde{P}_{\lambda_1\delta}\ast \tilde{P}_{\lambda_2\delta}\ast\cdots\ast \tilde{P}_{\lambda_t\delta}
\end{displaymath}
and
$$S_{\lambda}=\det(P_{(\lambda_k-k+m)\delta})_{1\leq k,m\leq t}.$$
The relation between $\tilde{P}_{\lambda}$ and $S_{\lambda}$ is  $$\tilde{P}_{\lambda}=\sum_{\mu\in P}K_{\lambda\mu}S_{\mu},$$
where $K_{\lambda\mu}$ is the Kostka number associated to the partitions $\lambda$ and $\mu$.

\begin{theorem}[\cite{Mcgerty2004The_Kronecker_quiver}]
The set
$\{S_{\lambda}\,\,|\,\,\textrm{$\lambda$ is a partition}\}$ is a subset of the canonical basis $\mathbf{B}$.
\end{theorem}

Let $\mathcal{G}$ be the set of $(c=(c_{-},c_{+}),\lambda)$, where $c_{-}:Prep\rightarrow\mathbb{N}$, $c_{+}:Prei\rightarrow\mathbb{N}$ are functions with finite support and $\lambda$ is a partition.

For any $(c,\lambda)\in\mathcal{G}$, consider
$$N'(c,\lambda)=\langle{M(c_{-})}\rangle\ast \tilde{P}_{\lambda}\ast\langle{M(c_{+})}\rangle$$
and
$$N(c,\lambda)=\langle{M(c_{-})}\rangle\ast S_{\lambda}\ast\langle{M(c_{+})}\rangle$$
where
\begin{displaymath}
M(c_{-})=\bigoplus_{\alpha\in Prep}M_\alpha^{\oplus c_{-}(\alpha)}
\end{displaymath}
and
\begin{displaymath}
M(c_{+})=\bigoplus_{\alpha\in Prei}M_\alpha^{\oplus c_{+}(\alpha)}.
\end{displaymath}

\begin{proposition}[\cite{Chen-Rootvector}\cite{Zhang2000PBW}]
The sets
$\{N(c,\lambda)\,\,|\,\,(c,\lambda)\in\mathcal{G}\}$ and $\{N'(c,\lambda)\,\,|\,\,(c,\lambda)\in\mathcal{G}\}$
are two $\mathbb{Z}[v,v^{-1}]$-bases of $\mathcal{C}^{\ast}(Q)_{\mathbb{Z}[v,v^{-1}]}$.
\end{proposition}

The sets
$\{N(c,\lambda)\,\,|\,\,(c,\lambda)\in\mathcal{G}\}$ and $\{N'(c,\lambda)\,\,|\,\,(c,\lambda)\in\mathcal{G}\}$
are called PBW-type bases of $\mathcal{C}^{\ast}(Q)$.

\section{The construction for cyclic quivers}\label{DDX}

The construction of various bases of affine $A$ type was obtained in \cite{DengDuXiao2007Generic} by considering the Hall algebra of the cyclic quiver. Let $Q$ be the following cyclic quiver whose vertex set is $I=\{0,1,2,\ldots,n\}$:
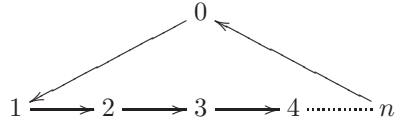
$$\xymatrix{
        &&0\ar[dll]&\\
1\ar[r]&2\ar[r]&3\ar[r]& 4\ar@{.}[r]&n\ar[ull]
}$$
Denote by $\mathcal{H}^{\ast}$ the twisted Ringel-Hall algebra of the category of nilpotent representations of $Q$ and $\mathcal{C}^{\ast}$ the twisted composition subalgebra of $\mathcal{H}^{\ast}$.

A multisegment is a formal sum
$$\pi=\sum_{i\in I,l\geq1}\pi_{il}[i,l],$$
where $\pi_{il}\in\mathbb{N}$ and $\{i\in I,l\geq1\,\,|\,\,\pi_{il}\neq0\}$ is a finite set. Let $\Pi$ be the set of multisegments.

There is a bijection between the set $\Pi$ and the isomorphism classes of nilpotent representations of $Q$. The isomorphism classes corresponding to $\pi$ is
$$
M(\pi)=\bigoplus_{i\in I,l\geq 1}S_i[l]^{\oplus \pi_{il}},
$$
where $S_i[l]$ is the unique indecomposable $kQ$-module with top $S_i$ and length $l$.

An element $\pi\in\Pi$ is called aperiodic, if $$\prod_{i\in I}\pi_{il}=0$$ for each $l\geq1$. The set of all aperiodic multisegments is denoted by $\Pi^{a}$.

There is a partial order on $\Pi$ defined as follows: for $\pi',\pi\in\Pi$ with the same dimension vector, $\pi'<\pi$ if and only if $\dim\textrm{Hom}(M,M(\pi'))>\dim\textrm{Hom}(M,M(\pi))$ for all
 indecomposable nilpotent representations $M$ of $Q$.

\begin{proposition}[\cite{DengDuXiao2007Generic}]
For each $\pi\in\Pi^{a}$, there exists a monomial $m^{\omega_{\pi}}$ on the divided powers of $u_i$ such that
$$m^{\omega_\pi}=\langle{M(\pi)}\rangle+\sum_{\pi'<\pi}\eta^{\pi'}_{\pi}\langle{M(\pi')}\rangle$$
with $\eta_{\pi}^{\pi'}\in\mathbb{Z}[v,v^{-1}]$.
\end{proposition}

Every non-empty subset of $\Pi^{a}$ contains a minimal element.
Define $E_{\pi}$ for all $\pi\in\Pi^a$ inductively by the following relations.
If $\pi\in\Pi^{a}$ is minimal,
\begin{displaymath}
E_{\pi}={m}^{w_{\pi}}=\langle{M(\pi)}\rangle+\sum_{\pi'<\pi,\pi'\in\Pi\backslash\Pi^a}\eta_{\pi}^{\pi'}\langle{M(\pi')}\rangle,
\end{displaymath}
If $E_{\pi'}$ have been defined for all $\pi>\pi'\in\Pi^a$, then
\begin{eqnarray*}
E_{\pi}&=&{m}^{w_{\pi}}-\sum_{\pi'<\pi,\pi'\in\Pi^a}\eta_{\pi}^{\pi'}E_{\pi'}\\
&=&
\langle{M(\pi)}\rangle+\sum_{\pi'<\pi,\pi'\in\Pi\backslash\Pi^a}\gamma_{\pi}^{\pi'}\langle{M(\pi')}\rangle.
\end{eqnarray*}

\begin{proposition}[\cite{DengDuXiao2007Generic}]
The set $\{{E}_\pi\,\,|\,\,\pi\in\Pi^a\}$ is a $\mathbb{Z}[v,v^{-1}]$ basis of $\mathcal{C}^{\ast}_{\mathbb{Z}[v,v^{-1}]}$, satisfying the following conditions
\begin{enumerate}
  \item[(1)]$\{{E}_\pi\,\,|\,\,\pi\in\Pi^a\}$ is independent of the choice of monomials ${m}^{w_{\pi}}$;
  \item[(2)]$$\overline{{E}_\pi}={E}_\pi+\sum_{\pi'<\pi}r^{\pi'}_{\pi}{E}_{\pi'}$$ with $r^{\pi'}_{\pi}\in\mathbb{Z}[v,v^{-1}]$.
\end{enumerate}
\end{proposition}

The set
$\{{E}_\pi\,\,|\,\,\pi\in\Pi^a\}$
is called a PBW-type basis of $\mathcal{C}^{\ast}$.

By Lemma \ref{Lemma_Lusztig}, there exists a unique family of elements $p^{\pi'}_{\pi}\in\mathbb{Z}[v^{-1}]$ defined for all $\pi'\leq\pi$ in $\Pi^a$ such that
\begin{enumerate}
  \item[(1)]$p^{\pi}_{\pi}=1$ for all $\pi\in\Pi^a$;
  \item[(2)]$p^{\pi'}_{\pi}\in v^{-1}\mathbb{Z}[v^{-1}]$ for all $\pi'\leq\pi$ in $\Pi^a$;
  \item[(3)]for all $\pi'\leq\pi$ in $\Pi^a$,
  $$p^{\pi'}_{\pi}=\sum_{\pi'',\pi'\leq \pi''\leq \pi}\overline{p^{\pi'}_{\pi''}}r^{\pi''}_{\pi}.$$
\end{enumerate}
For any $\pi\in\Pi^a$, let $$c_{\pi}=E_{\pi}+\sum_{\pi'<\pi}p^{\pi'}_{\pi}{E}_{\pi'}.$$

\begin{theorem}[\cite{DengDuXiao2007Generic}]
The set $\{c_{\pi}\,\,|\,\,\pi\in\Pi^a\}$ is a $\mathbb{Z}[v,v^{-1}]$-basis of $\mathcal{C}^{\ast}_{\mathbb{Z}[v,v^{-1}]}$
satisfying the following conditions.
\begin{enumerate}
  \item[(1)]$\overline{c_{\pi}}=c_{\pi}$;
  \item[(2)]$(c_{\pi},c_{\pi'})\in\delta_{\pi,\pi'}+v^{-1}\mathbb{Z}[[v^{-1}]]\cap\mathbb{Q}(v)$.
\end{enumerate}
\end{theorem}

\begin{corollary}[\cite{DengDuXiao2007Generic}]
The set $\{c_{\pi}\,\,|\,\,\pi\in\Pi^a\}$ is the canonical basis of $\mathcal{C}^{\ast}$.
\end{corollary}

\section{The construction for tame quivers I}\label{LXZ}

This construction of bases of affine $A,D,E$ type was obtained in \cite{Lin_Xiao_Zhang_Representations_of_tame_quivers_and_affine_canonical_bases} by using the Ringel-Hall algebra approach. Let $Q$ be an acyclic quiver of affine type.
Give an order on $I=\{0,1,2,\ldots,n\}$ such that $i>j$ implies that there doesn't exist an arrow $i\rightarrow j$. Define a double infinite sequence $$h=(\cdots,i_{-1},i_0,i_1,\cdots)$$
such that $i_k=k$ for all $k=0,1,2,\ldots,n$ and $i_{k+n+1}=i_k$ for all $k\in\mathbb{Z}$.
Then $$\{\beta_{0}=\alpha_{i_0},\beta_{-1}=s_{i_0}(\alpha_{i_{-1}}),\beta_{-2}=s_{i_0}s_{i_{-1}}(\alpha_{i_{-2}}),\cdots\}$$ is the set of dimension vectors of all indecomposable preprojective modules and
$$\{\beta_{1}=\alpha_{i_1},\beta_{2}=s_{i_1}(\alpha_{i_2}),\beta_{3}=s_{i_1}s_{i_2}(\alpha_{i_3}),\cdots\}$$
is the set of dimension vectors of all indecomposable preinjective modules.

The category $\textrm{rep}_kQ$ has direct sum decomposition
$$\textrm{rep}_kQ=Prep\oplus Reg\oplus Prei$$ and each component is closed on taking extensions
in $\textrm{rep}_kQ$ and direct summands.
Thus the Hall algebras of these components are
subalgebras of $\mathcal{H}^{\ast}(kQ)$. Generic composition algebra $\mathcal{C}^{\ast}(Q)$ contains the Hall algebras of the components
$Prep$ and $Prei$ as subalgebras. Under the isomorphism between $\mathcal{C}^{\ast}(Q)$ and $\mathbf{U}^+$, we can view them as subalgebras of $\mathbf{U}^+$.

Let
$$\langle{M(\beta_k)}\rangle=\left\{\begin{array}{cc}
                             T^{-1}_{i_0}T^{-1}_{i_{-1}}\cdots T^{-1}_{i_{k+1}}(E_{i_k}) & \textrm{if $k\leq0$,}\\
                             T_{1}T_{i_{2}}\cdots T_{i_{k-1}}(E_{i_k}) & \textrm{if $k>0$.}\\
                           \end{array}
\right.$$
For a support finite function $a:\mathbb{Z}^{\leq0}\rightarrow\mathbb{N}$, define
\begin{eqnarray*}
\langle{M(a)}\rangle&=&\langle{\oplus_{k\leq0}M(\beta_k)^{\oplus{a}(k)}}\rangle\\
&=&E_{i_0}^{(a(0))}T^{-1}_{i_0}(E_{i_{-1}}^{(a(-1))})T^{-1}_{i_0}T^{-1}_{i_{-1}}(E_{i_{-2}}^{(a(-2))})\cdots.
\end{eqnarray*}
For a support finite function $b:\mathbb{Z}^{>0}\rightarrow\mathbb{N}$, define
\begin{eqnarray*}\langle{M(b)}\rangle&=&\langle{\oplus_{k>0}M(\beta_k)^{\oplus{b}(k)}}\rangle\\&=&\cdots T_{i_1}T_{i_2}(E_{i_3}^{(b(3))})T_{i_1}(E_{i_2}^{(b(2))})E_{i_1}^{(b(1))}.\end{eqnarray*}
Note that $\langle{M(a)}\rangle$ and $\langle{M(b)}\rangle$ belong to the Hall algebras of the components
$Prep$ and $Prei$ respectively.

For regular part, there exist $s (s\leq3)$ non-homogeneous  tubes $J_1,J_2,\ldots,J_s$.
It is well-known that the full subcategory corresponding to the tube $J_t$ is equivalent to the category of nilpotent representations
of the cyclic quiver with $r_t$ vertices, where $r_t$ is the rank of $J_t$.
Choose such an equivalence for each $J_t$, which induces an algebra isomorphism
$$\varepsilon_t:\mathcal{H}_{r_t}^{\ast}\rightarrow\mathcal{H}^{\ast}(J_t),$$
where $\mathcal{H}_{r_t}^{\ast}$ is the twisted Ringel-Hall algebra corresponding to the cyclic quiver with $r_t$ vertices
and $\mathcal{H}^{\ast}(J_t)$ is the twisted Ringel-Hall algebra of the full subcategory corresponding to the tube $J_t$.

In Section \ref{DDX}, a PBW-type basis $\{E_{\pi}\,\,|\,\,\pi\in\Pi^{a}_{t}\}$ for the composition subalgebra $\mathcal{C}_{r_t}^{\ast}$ of the twisted Ringel-Hall algebra $\mathcal{H}_{r_t}^{\ast}$ has been constructed.
For any $\pi\in\Pi^a_{t}$, the image of $E_{\pi}$ under $\varepsilon_t$ is still denoted by
$E_{\pi}$.

Let $K_2$ be the path algebra of the Kronecker quiver
and  $F:\textrm{mod}kK_2\hookrightarrow\textrm{mod}kQ$ be the canonical embedding.
This embedding induces a map $F:\mathcal{H}^{\ast}(kK_2)\rightarrow\mathcal{H}^{\ast}(kQ)$.

In Section \ref{Kronecker}, $\tilde{\Psi}_{n\delta}$ and $\tilde{P}_{n\delta}$ have been defined. Let
$$E_{n\delta}=F(\tilde{P}_{n\delta})$$ for $n\in\mathbb{Z}_{>0}$ and $$E_{p\delta}=E_{p_1\delta}\ast\cdots\ast E_{p_s\delta}$$ for a partition $p=(p_1\geq \cdots\geq p_s)$.

Let $\mathcal{M}$ be the set of quadruples $c=(a_c,b_c,\pi_c,p_c)$, where $a_c:\mathbb{Z}^{\leq0}\rightarrow\mathbb{N}$ and $b_c:\mathbb{Z}^{>0}\rightarrow\mathbb{N}$ are functions with finite support, $\pi_{c}\in\Pi^{a}_{1}\times\cdots\times\Pi^{a}_{s}$ and $p_c$ is a partition.
For each $c\in\mathcal{M}$, define
$$E^{c}=\langle{M(a_c)}\rangle\ast E_{\pi_{1}c}\ast\cdots\ast E_{\pi_{s}c}\ast E_{p_c\delta}\ast\langle{M(b_c)}\rangle.$$

Recall that $E_\nu$ is the variety consisting of representations with dimension vector $\nu\in\mathbb{N}I$ of the quiver $Q$ over $\bar{k}$. For subset $\mathcal{A}\subset E_\alpha$ and $\mathcal{B}\subset E_\beta$, define the extension set $\mathcal{A}\star\mathcal{B}$ of $\mathcal{A}$ by $\mathcal{B}$ to be
the set of all $z\in E_{\alpha+\beta}$ such that $M(z)$ is an extension of $M(x)$ by $M(y)$ for some $x\in\mathcal{A}, y\in\mathcal{B}$.

Define the subvariety of $E_\nu$
$$\mathcal{O}_c=\mathcal{O}_{M(a_c)}\star\mathcal{O}_{M(\pi_{1}c)}\star\mathcal{O}_{M(\pi_{2}c)}\star\dots\star\mathcal{O}_{M(\pi_{s}c)}\star\mathcal{N}_{p_c}\star\mathcal{O}_{M(b_c)}$$
for any $c\in\mathcal{M}$, where $\mathcal{N}_{p}=\mathcal{N}_{p_1}\star\mathcal{N}_{p_2}\star\dots\star\mathcal{N}_{p_s}$ if $p=(p_1\geq \cdots\geq p_s)$ and
$\mathcal{N}_{p_i}$ are the union of orbits of images of all regular modules in $kK_2$ under $F$ with dimension vector $p_i\delta$.

\begin{proposition}[\cite{Lin_Xiao_Zhang_Representations_of_tame_quivers_and_affine_canonical_bases}]
The set $\{E^c\,\,|\,\,c\in\mathcal{M}\}$ is a $\mathbb{Q}(v)$-basis of $\mathcal{C}^{\ast}(kQ)$.
\end{proposition}

The set $\{E^c\,\,|\,\,c\in\mathcal{M}\}$ is a PBW-type basis of $\mathcal{C}^{\ast}(kQ)$.

\begin{proposition}[\cite{Lin_Xiao_Zhang_Representations_of_tame_quivers_and_affine_canonical_bases}]
For each $c\in\mathcal{M}$, there exists a monomial $m_c$ on the divided powers of $u_i$ such that
$$m_c=E^c+\sum_{\dim\mathcal{O}_{c'}<\dim\mathcal{O}_{c}}h_{c}^{c'}E^{c'}$$
with $h_{c}^{c'}\in\mathbb{Q}[v,v^{-1}]$.
\end{proposition}

Similarly to the case of finite type, from this basis we can get a bar-invariant basis. But it is not the canonical basis considered by Lusztig. Hence in \cite{Lin_Xiao_Zhang_Representations_of_tame_quivers_and_affine_canonical_bases}, another PBW-type basis is constructed.

There is a bilinear form $(-,-)$ on $\mathcal{H}_q^{\ast}(kQ)$ defined in \cite{green1995hall}.
Consider the $\mathbb{Q}(v)$-basis $\{E^{{c}}\,\,|\,\,{c}\in\mathcal{M}\}$.
Let $R(\mathcal{C}^{\ast}(kQ))$ be the $\mathbb{Q}(v)$-subspace of $\mathcal{C}^{\ast}(kQ)$ with the basis $\{E_{\pi_{1\mathbf{c}}}\ast E_{\pi_{2{c}}}\ast\cdots\ast E_{\pi_{s{c}}}\ast E_{p_{{c}}\delta}\}$,
where $\pi_{{c}}=(\pi_{1{c}},\pi_{2{c}},\ldots,\pi_{s{c}})\in\Pi_{1}^a\times\Pi_{2}^a\times\cdots\times\Pi_{s}^a$, and $p_{{c}}=({p_1}\geq{p_2}\geq\cdots\geq{p_t})$ is a partition. Note that $R(\mathcal{C}^{\ast}(kQ))$ is a subalgebra of $\mathcal{C}^{\ast}(kQ)$.

Let $R^a(\mathcal{C}^{\ast}(kQ))$ be the subalgebra of $R(\mathcal{C}^{\ast}(kQ))$ with the basis $\{E_{\pi_{1{c}}}\ast E_{\pi_{2{c}}}\ast\cdots\ast E_{\pi_{s{c}}}\,\,|\,\,\pi_{{c}}=(\pi_{1{c}},\pi_{2{c}},\ldots,\pi_{s{c}})\in\Pi_{1}^a\times\Pi_{2}^a\times\cdots\times\Pi_{s}^a\}$. For any $\alpha,\beta\in\mathbb{N}I$, define $\alpha\leq\beta$ if $\beta-\alpha\in\mathbb{N}I$. If $\beta<\delta$, $R(\mathcal{C}^{\ast}(kQ))_{\beta}=R^a(\mathcal{C}^{\ast}(kQ))_{\beta}$. Define $\mathcal{F}_{\delta}=\{x\in R(\mathcal{C}^{\ast}( kQ))_{\delta}\,\,|\,\,(x,R^a(\mathcal{C}^{\ast}(kQ))_{\delta})=0\}$.

In \cite{Lin_Xiao_Zhang_Representations_of_tame_quivers_and_affine_canonical_bases}, it is proved that
\begin{displaymath}
R(\mathcal{C}^{\ast}(kQ))_{\delta}=R^a(\mathcal{C}^{\ast}(kQ))_{\delta}\oplus\mathcal{F}_{\delta}
\end{displaymath}
and $\dim\mathcal{F}_{\delta}=1$. By the method of Schmidt orthogonalization, we may set
\begin{displaymath}
E'_{\delta}=E_{\delta}-\sum_{M(\pi_{i\mathbf{c}}),\underline{\dim} M(\pi_{i{c}})=\delta, 1\leq i\leq s}a_{\pi_{i{c}}}E_{\pi_{i{c}}}.
\end{displaymath}
Then $\mathcal{F}_{\delta}=\mathbb{Q}(v)E'_{\delta}$.

Now let $R(\mathcal{C}^{\ast}(kQ))(1)$ be the subalgebra of $R(\mathcal{C}^{\ast}(kQ))$ generated by $R^a(\mathcal{C}^{\ast}(kQ))$ and $\mathcal{F}_{\delta}$. If $\beta<2\delta$, $R(\mathcal{C}^{\ast}(kQ))(1)_{\beta}=R(\mathcal{C}^{\ast}(kQ))_{\beta}$. Define
\begin{displaymath}
\mathcal{F}_{2\delta}=\{x\in R(\mathcal{C}^{\ast}(kQ))_{2\delta}\,\,|\,\,(x,R(\mathcal{C}^{\ast}(kQ))(1)_{2\delta})=0\}.
\end{displaymath}
Then $\dim\mathcal{F}_{2\delta}=1$ and $R(\mathcal{C}^{\ast}(kQ))_{2\delta}=R(\mathcal{C}^{\ast}(kQ))(1)_{2\delta}\oplus\mathcal{F}_{2\delta}$.

Suppose $R(\mathcal{C}^{\ast}(kQ))(n-1)$ has been defined, we define
\begin{displaymath}
\mathcal{F}_{n\delta}=\{x\in R(\mathcal{C}^{\ast}(kQ))_{n\delta}\,\,|\,\,(x,R(\mathcal{C}^{\ast}(kQ))(n-1)_{n\delta})=0\}.
\end{displaymath}
Let $R(\mathcal{C}^{\ast}(kQ))(n)$ be the subalgebra of $R(\mathcal{C}^{\ast}(kQ))$ generated by $R(\mathcal{C}^{\ast}(kQ))(n-1)$ and $\mathcal{F}_{n\delta}$. Then $\dim\mathcal{F}_{n\delta}=1$ and $R(\mathcal{C}^{\ast}(kQ))_{n\delta}=R(\mathcal{C}^{\ast}(kQ))(n-1)_{n\delta}\oplus\mathcal{F}_{n\delta}$. Similarly, choose $E'_{n\delta}$ such that $E_{n\delta}-E'_{n\delta}\in R(\mathcal{C}^{\ast}(kQ))(n-1)_{n\delta}$ and $\mathcal{F}_{n\delta}=\mathbb{Q}(v)E'_{n\delta}$ for all $n\geq1$.

Let $P_{n\delta}=nE'_{n\delta}$ and $$S'_{p_c\delta}=\det(P_{((p_c)_k-k+m)\delta})_{1\leq k,m\leq t}$$
be the Schur functions corresponding to $P_{n\delta}$.

For each $c\in\mathcal{M}$, define
$$e^{c}=\langle{M(a_c)}\rangle\ast E_{\pi_{1}c}\ast\cdots\ast E_{\pi_{s}c}\ast S'_{p_c\delta}\ast\langle{M(b_c)}\rangle.$$
The set $\{e^{c}\,\,|\,\,c\in\mathcal{M}\}$ is another PBW-type basis of $\mathcal{C}^{\ast}(kQ)$.

For two $c,c'\in\mathcal{M}$, define $e^{c'}<e^{c}$ if either $\dim\mathcal{O}_{c'}<\dim\mathcal{O}_{c}$ or $\dim\mathcal{O}_{c'}=\dim\mathcal{O}_{c}$ but $p_c<p_{c'}$ under lexicographic order of partitions.

\begin{proposition}[\cite{Lin_Xiao_Zhang_Representations_of_tame_quivers_and_affine_canonical_bases}]
The set
$\{e^c\,\,|\,\,c\in\mathcal{M}\}$ is a $\mathbb{Q}[v,v^{-1}]$-basis of $\mathcal{C}^{\ast}(kQ)_{\mathbb{Q}[v,v^{-1}]}$ satisfying
\begin{enumerate}
  \item[(1)]$(e^c,e^{c'})\in\delta_{c,c'}+v^{-1}\mathbb{Q}[[v^{-1}]]\cap\mathbb{Q}(v)$;
  \item[(2)]$$m_c=e^c+\sum_{e^{c'}<e^{c}}a_{c}^{c'}e^{c'}$$with $a_{c}^{c'}\in\mathbb{Q}[v,v^{-1}]$.
\end{enumerate}
\end{proposition}
Similarly to the case of finite type, Lin-Xiao-Zhang proved the following Theorem  by using Lemma \ref{Lemma_Lusztig}.
\begin{theorem}[\cite{Lin_Xiao_Zhang_Representations_of_tame_quivers_and_affine_canonical_bases}]
For any $c\in\mathcal{M}$, there exists a unique $\mathcal{E}^c\in\mathcal{C}^{\ast}(Q)_{\mathbb{Q}[v,v^{-1}]}$ satisfying the following conditions
\begin{enumerate}
  \item[(1)]$\overline{\mathcal{E}^c}=\mathcal{E}^c$;
  \item[(2)]$(\mathcal{E}^c,\mathcal{E}^{c'})\in\delta_{c,c'}+v^{-1}\mathbb{Q}[[v^{-1}]]\cap\mathbb{Q}(v)$;
  \item[(3)]$$\mathcal{E}^c=e^c+\sum_{e^{c'}<e^{c}}b_{c}^{c'}e^{c'}$$ with $b_{c}^{c'}\in\mathbb{Q}[v,v^{-1}]$.
\end{enumerate}
Moreover, the set $\{\mathcal{E}^c\,\,|\,\,c\in\mathcal{M}\}$ is  the canonical basis of $\mathcal{C}^{\ast}(Q)$.
\end{theorem}

\section{The construction for tame quivers II}

  Let $Q$ be an  acyclic quiver of affine type.
Denoted by $J_1,J_2,\ldots,J_s(s\leq3)$ the non-homogeneous tubes.
Let $\mathcal{H}^0(Q)$ be the $\mathbb{Q}(v)$-subalgebra of $\mathcal{H}^{\ast}(Q)$ generated by $u_i$ for $i\in I$ and $u_{[M]}$ for $M\in J_i$. Note that  $\mathcal{C}^{\ast}(Q)\subset\mathcal{H}^0(Q)$ and $\mathcal{H}^0(Q)$ is called the reductive extension of $\mathcal{C}^{\ast}(Q)$.

With the same notations in Section \ref{LXZ}, there is a double infinite sequence $$h=(\cdots,i_{-1},i_0,i_1,\cdots)$$ such that
$$\{\beta_{0}=\alpha_{i_0},\beta_{-1}=s_{i_0}(\alpha_{i_{-1}}),\beta_{-2}=s_{i_0}s_{i_{-1}}(\alpha_{i_{-2}}),\cdots\}$$ is the set of dimension vectors of all indecomposable preprojective modules and
$$\{\beta_{1}=\alpha_{i_1},\beta_{2}=s_{i_1}(\alpha_{i_2}),\beta_{3}=s_{i_1}s_{i_2}(\alpha_{i_3}),\cdots\}$$
is the set of dimension vectors of all indecomposable preinjective modules.
We order these $\beta_{t}$ for various $t\in\mathbb{Z}$ by
$$(\beta_0<\beta_{-1}<\beta_{-2}<\cdots)<(\cdots\beta_3<\beta_2<\beta_1).$$

For a support finite function $c_{-}:\mathbb{Z}^{\leq0}\rightarrow\mathbb{N}$, define
\begin{eqnarray*}
\langle{M(c_{-})}\rangle&=&\langle{\oplus_{k\leq0}M(\beta_k)^{\oplus{c_{-}}(k)}}\rangle\\
&=&E_{i_0}^{(c_{-}(0))}T^{-1}_{i_0}(E_{i_{-1}}^{(c_{-}(-1))})T^{-1}_{i_0}T^{-1}_{i_{-1}}(E_{i_{-2}}^{(c_{-}(-2))})\cdots.
\end{eqnarray*}
For a support finite function $c_{+}:\mathbb{Z}^{>0}\rightarrow\mathbb{N}$, define
\begin{eqnarray*}
\langle{M(c_{+})}\rangle&=&\langle{\oplus_{k>0}M(\beta_k)^{\oplus{c_{+}}(k)}}\rangle\\
&=&\cdots T_{i_1}T_{i_2}(E_{i_3}^{(c_{+}(3))})T_{i_1}(E_{i_2}^{(c_{+}(2))})E_{i_1}^{(c_{+}(1))}.
\end{eqnarray*}

For any $c_{0}=(\pi_1,\ldots,\pi_s)\in\Pi_{1}\times\cdots\times\Pi_{s}$, let
$$M(c_{0})=\varepsilon_1(M(\pi_1))\ast\cdots\ast\varepsilon_s(M(\pi_s)),$$
where
$$\varepsilon_t:\mathcal{H}_{r_t}^{\ast}\rightarrow\mathcal{H}^{\ast}(J_t)$$
is the algebra isomorphism in Section \ref{LXZ}.

Similarly to the case of Kronecker quiver, define
$$\tilde{P}_{m}=\sum_{[M]:\underline{\dim}M=m\delta\atop\textrm{$M$ is homogeneous regular}}v^{-\dim M}u_{[M]}$$
for $m\in\mathbb{Z}_{\geq0}$ and
$$\tilde{P}_{\lambda}=\prod_{1\leq k\leq t}\tilde{P}_{\lambda_k},$$
$$S_{\lambda}=\det(P_{\lambda_k-k+m})_{1\leq k,m\leq t}$$
for a partition $\lambda=(\lambda_1\geq\lambda_2\geq\cdots\geq\lambda_t)$.

We can write in terms of modules
 $$\tilde{P}_{\lambda}=\sum_{[M]:\underline{\dim}M=m\delta\atop\textrm{$M$ is homogeneous regular}}A_\lambda^{[M]}(v)\langle{M}\rangle,$$
  $$S_{\lambda}=\sum_{[M]:\underline{\dim}M=m\delta\atop\textrm{$M$ is homogeneous regular}}B_\lambda^{[M]}(v)\langle{M}\rangle,$$
where $A_\lambda^{[M]}(v),B_\lambda^{[M]}(v)\in\mathbb{Z}[v,v^{-1}]$.

It is interesting to compute $A_\lambda^{[M]}(v),B_\lambda^{[M]}(v)$ for some special homogeneous regular $M$.

Let $\mathcal{Z}_k$ be the set of all homogeneous tubes of $\textrm{mod}\, kQ$ and by $\deg z$ we denote the degree of the corresponding irreducible polynomial of $z\in\mathcal{Z}_k$. We denote by $M(l,z)$ the indecomposable module in tube $z$ with dimension vector $ld_z\delta$ where $d_z=\deg z$, that is, $l$ means the "level" of the corresponding module.

For a partition $\mu=(\mu_1\geq\mu_2\geq\cdots\geq\mu_t)$ and $\underline{z}=(z_1,z_2,\cdots,z_t)$ such that $z_i\in\mathcal{Z}_k$ and $\deg z_i=1$ for all $i$,  we denote
$$M(\mu,\underline{z})=M(\mu_1,z_1)\oplus M(\mu_2,z_2)\oplus\cdots\oplus M(\mu_t,z_t).$$

For a partition $\mu=(\mu_1\geq\mu_2\geq\cdots\geq\mu_t)$ and $\underline{z}'=(z'_1,z'_2,\cdots,z'_t)$ such that $z'_i\in\mathcal{Z}_k$ and $\deg z'_i=\mu_i$ for all $i$, we denote
$$M[\mu,\underline{z}']=M(1,z'_1)\oplus M(1,z'_2)\oplus\cdots\oplus M(1,z'_t).$$

Note that both $M(\mu,\underline{z})$ and $M[\mu,\underline{z}']$ have the dimension vector $|\mu|\delta$.

Let $K_{\lambda\mu}\in\mathbb{Z}$ be the Kostka numbers.

\begin{proposition}[\cite{Xiao-Xu-Zhao}]\label{Kostka}
  For partitions $\lambda,\mu$ with $|\lambda|=|\mu|$ and $\underline{z}=(z_1,z_2,\cdots,z_t)$ such that $z_i\in\mathcal{Z}_k$ and $\deg z_i=1$ for all $i$, $B_\lambda^{[M(\mu,z)]}(v)=v^{-|\lambda|}K_{\mu\lambda}$.
\end{proposition}

For a partition $\lambda$ of $m$, let $S^\lambda$ be the Specht
module for $\mathcal{S}_m$. Let $t_\lambda(\mu)=\chi_{S^\lambda}(g_\mu)$ be the complex character value of $S^\lambda$ at $g_\mu\in\mathcal{S}_m$
of cycle type $\mu$.
Then $(t_\lambda(\mu))_{\lambda,\mu}$ is the character table of $\mathcal{S}_m$. Let $t'_\lambda$ be the character of the permutation module $M^\lambda$. It is known that $t'_\lambda=\sum_{\mu}K_{\lambda\mu}t_\mu$.

\begin{proposition}[\cite{Xiao-Xu-Zhao}]
  For partitions $\lambda,\mu$ with $|\lambda|=|\mu|$ and $\underline{z}'=(z'_1,z'_2,\cdots,z'_t)$ such that $z'_i\in\mathcal{Z}_k$ and $\deg z'_i=\mu_i$ for all $i$, $A_\lambda^{[M[\mu,z']]}(v)=v^{-|\lambda|}t'_\lambda(\mu)$.
\end{proposition}
\begin{corollary}[\cite{Xiao-Xu-Zhao}]
  For partitions $\lambda,\mu$ with $|\lambda|=|\mu|$ and $\underline{z}'=(z'_1,z'_2,\cdots,z'_t)$ such that $z'_i\in\mathcal{Z}_k$ and $\deg z'_i=\mu_i$ for all $i$, $B_\lambda^{[M[\mu,z']]}(v)=v^{-|\lambda|}t_\lambda(\mu)$.
\end{corollary}

Let $\mathcal{G}$ be the set of $(c=(c_{-},c_{0},c_{+}),t_\lambda)$, where $c_{-}:\mathbb{Z}^{\leq0}\rightarrow\mathbb{N}$ (resp. $c_{+}:\mathbb{Z}^{>0}\rightarrow\mathbb{N}$) is function with finite support, $c_{0}\in\Pi_{1}\times\cdots\times\Pi_{s}$ and $t_\lambda$ is the character of a Specht module $S^\lambda$.
Let $\mathcal{G}^a$ the subset of $\mathcal{G}$ consisting of all such $(c,t_\lambda)$ such that $c_{0}\in\Pi^a_{1}\times\cdots\times\Pi^a_{s}$.

For any $(c,t_\lambda)\in\mathcal{G}$, consider
$$N(c,t_\lambda)=\langle{M(c_{-})}\rangle\ast \langle{M(c_{0})}\rangle\ast S_{\lambda}\ast\langle{M(c_{+})}\rangle,$$
where $c=(c_{-},c_{0},c_{+})$.

\begin{proposition}[\cite{Xiao-Xu-Zhao}]
The set $\{N(c,t_\lambda)\,\,|\,\,(c,t_\lambda)\in\mathcal{G}\}$ is an $\mathbb{Q}(v)$-basis of $\mathcal{H}^0$ such that
\begin{enumerate}
  \item[(1)]$(N(c,t_\lambda),N(c',t_\lambda'))\in\delta_{(c,t_\lambda),(c',t_\lambda')}+v^{-1}\mathbb{Q}[[v^{-1}]]\cap\mathbb{Q}(v)$;
  \item[(2)]$$N(c,t_\lambda)\ast N(c',t_\lambda')=\sum_{(c'',t_\lambda'')\in\mathcal{G}}P_{(c,t_\lambda),(c',t_\lambda')}^{(c'',t_\lambda'')}N(c'',t_\lambda'')$$ with $P_{(c,t_\lambda),(c',t_\lambda')}^{(c'',t_\lambda'')}\in\mathbb{Z}[v,v^{-1}]$.
\end{enumerate}
\end{proposition}

There is a "combinatorial" order $<$ on $\mathcal{G}$ defined as follows.
For $c_-,c_-':\mathbb{Z}^{\leq0}\rightarrow\mathbb{N}$, define $c_-<c_-'$ if and only if there exists $j\leq0$ such that $c_-(t)=c_-'(t)$ for all $j<t\leq0$ and $c_-(j)>c_-'(j)$. For $c_+,c_+':\mathbb{Z}^{>0}\rightarrow\mathbb{N}$, define $c_+<c_+'$ if and only if there exists $j>0$ such that $c_-(t)=c_-'(t)$ for all $j>t>0$ and $c_-(j)>c_-'(j)$.
The partial order on $\Pi$ is given in Section \ref{DDX}.
For $t_\lambda$ and $t_\lambda'$, $t_\lambda<t_{\lambda'}$ means that $\lambda$ is less than $\lambda'$ under lexicographic order of partitions.

\begin{definition}
For $(c,t_\lambda),(c',t_{\lambda'})\in\mathcal{G}$, let $c=(c_{-},c_{0},c_{+})$, $c'=(c'_{-},c'_{0},c'_{+})$, $c_{0}=(\pi_1,\ldots,\pi_s)$ and $c'_{0}=(\pi'_1,\ldots,\pi'_s)$.
Define $(c',t_{\lambda'})<(c,t_\lambda)$ if one of the following three conditions holds.
\begin{enumerate}
  \item[(1)]$c=c'$ and $t_{\lambda'}>t_\lambda$;
  \item[(2)]$c'_-\leq c_-$, $c'_+\leq c_+$ but not all equalities hold;
  \item[(3)]$c_-=c_-'$, $c_+=c_+'$, $\pi'_1\leq\pi_1,\ldots,\pi'_s\leq\pi_s$ but not all equalities hold.
\end{enumerate}
\end{definition}

\begin{proposition}[\cite{Xiao-Xu-Zhao}]
For each $(c,t_\lambda)\in\mathcal{G}^a$, there exists a monomial $m^{\omega(c,t_\lambda)}$ on the divided powers of $u_i$ such that
$$m^{\omega(c,t_\lambda)}=N(c,t_\lambda)+\sum_{(c',t_{\lambda'})<(c,t_\lambda)}a_{c,t_\lambda}^{c',t_{\lambda'}}N(c',t_{\lambda'})$$
with $a_{c,t_\lambda}^{c',t_{\lambda'}}\in\mathbb{Z}[v,v^{-1}]$.\end{proposition}

Li gave the geometric construction of this monomial basis in \cite{Li-yiqiang_canonical_monomial}.

With this partial order on $\mathcal{G}^a$, every nonempty subset has a minimal element.
Define $E{(c,t_\lambda)}$ for all $(c,t_\lambda)\in\mathcal{G}^a$ inductively by the following relations.
If $(c,t_\lambda)\in\mathcal{G}^a$ is minimal,
\begin{displaymath}
E{(c,t_\lambda)}=m^{\omega(c,t_\lambda)}=N(c,t_\lambda)+\sum_{(c',t_{\lambda'})<(c,t_\lambda)\atop(c',t_{\lambda'})\in\mathcal{G}\backslash\mathcal{G}^a}a_{c,t_\lambda}^{c',t_{\lambda'}}N(c',t_{\lambda'}).
\end{displaymath}
If $E{(c',t_{\lambda'})}$ have been defined for all $(c,t_\lambda)>(c',t_{\lambda'})\in\mathcal{G}^a$, then
\begin{eqnarray*}
E{(c,t_\lambda)}&=&m^{\omega(c,t_\lambda)}-\sum_{(c',t_{\lambda'})<(c,t_\lambda)\atop(c',t_{\lambda'})\in\mathcal{G}^a}a_{c,t_\lambda}^{c',t_{\lambda'}}E(c',t_{\lambda'})\\
&=&N(c,t_\lambda)+\sum_{(c',t_{\lambda'})<(c,t_\lambda)\atop(c',t_{\lambda'})\in\mathcal{G}\backslash\mathcal{G}^a}b_{c,t_\lambda}^{c',t_{\lambda'}}N(c',t_{\lambda'}).
\end{eqnarray*}

\begin{proposition}[\cite{Xiao-Xu-Zhao}]The set $\{E{(c,t_\lambda)}\,\,|\,\,(c,t_\lambda)\in\mathcal{G}^a\}$ is a $\mathbb{Z}[v,v^{-1}]$ basis of $\mathcal{C}^{\ast}(Q)$, such that
\begin{enumerate}
  \item[(1)]$\{E{(c,t_\lambda)}\,\,|\,\,(c,t_\lambda)\in\mathcal{G}^a\}$ is independent of the choice of monomials $m^{\omega(c,t_\lambda)}$.
  \item[(2)]$$\overline{E{(c,t_\lambda)}}=E{(c,t_\lambda)}+\sum_{(c',t_{\lambda'})<(c,t_\lambda)\atop(c',t_{\lambda'})\in\mathcal{G}^a}\gamma_{c,t_\lambda}^{c',t_{\lambda'}}E(c',t_{\lambda'})$$ with $\gamma_{c,t_\lambda}^{c',t_{\lambda'}}\in\mathbb{Z}[v,v^{-1}]$.
\end{enumerate}
\end{proposition}

The set $\{E{(c,t_\lambda)}\,\,|\,\,(c,t_\lambda)\in\mathcal{G}^a\}$ is called a PBW-type basis of $\mathcal{C}^{\ast}(Q)$.

By Lemma \ref{Lemma_Lusztig}, there exists a unique family of elements $\zeta_{c,t_\lambda}^{c',t_{\lambda'}}\in\mathbb{Z}[v^{-1}]$ defined for all $(c',t_{\lambda'})\leq (c,t_\lambda)$ in $\mathcal{G}^a$ such that
\begin{enumerate}
  \item[(1)]$\zeta_{c,t_\lambda}^{c,t_\lambda}=1$ for all $(c,t_\lambda)\in\mathcal{G}^a$;
  \item[(2)]$\zeta_{c,t_\lambda}^{c',t_{\lambda'}}\in v^{-1}\mathbb{Z}[v^{-1}]$ for all $(c',t_{\lambda'})\leq (c,t_\lambda)$ in $\mathcal{G}^a$;
  \item[(3)]for all $(c',t_{\lambda'})\leq (c,t_\lambda)$ in $\mathcal{G}^a$,$$\zeta_{c,t_\lambda}^{c',t_{\lambda'}}=\sum_{(c'',t_{\lambda''}),(c',t_{\lambda'})\leq (c'',t_{\lambda''})\leq (c,t_\lambda)}\overline{\zeta_{c',t_{\lambda'}}^{c'',t_{\lambda''}}}\gamma_{c,t_\lambda}^{c'',t_{\lambda''}}.$$
\end{enumerate}
For any $(c,t_\lambda)\in\mathcal{G}^a$, let $$C{(c,t_\lambda)}=E{(c,t_\lambda)}+\sum_{(c',t_{\lambda'})<(c,t_\lambda)\atop(c',t_{\lambda'})\in\mathcal{G}^a}\zeta_{c,t_\lambda}^{c',t_{\lambda'}}E{(c',t_{\lambda'})}.$$

\begin{theorem}[\cite{Xiao-Xu-Zhao}]
The set $\{C{(c,t_\lambda)}\,\,|\,\,(c,t_\lambda)\in\mathcal{G}^a\}$ is a $\mathbb{Z}[v,v^{-1}]$-basis of $\mathcal{C}^{\ast}(Q)_{\mathbb{Z}[v,v^{-1}]}$
satisfying the following conditions.
\begin{enumerate}
  \item[(1)]$\overline{C{(c,t_\lambda)}}=C{(c,t_\lambda)}$;
  \item[(2)]$(C{(c,t_\lambda)},C{(c',t_{\lambda'})})\in\delta_{(c,t_\lambda),(c',t_{\lambda'})}+v^{-1}\mathbb{Z}[[v^{-1}]]\cap\mathbb{Q}(v)$.
\end{enumerate}
\end{theorem}

\begin{corollary}[\cite{Xiao-Xu-Zhao}]
The set $\{C{(c,t_\lambda)}\,\,|\,\,(c,t_\lambda)\in\mathcal{G}^a\}$ is the canonical basis of $\mathcal{C}^{\ast}(Q)$.
\end{corollary}

\bibliography{mybibfile}

\begin{thebibliography}{10}

\bibitem{BCP}
J.~Beck, V.~Chari, and A.~Pressley.
\newblock An algebraic characterization of the affine canonical basis.
\newblock {\em Duke Math. J.}, 99(3):455--487, 1998.

\bibitem{2002Beck-Nakajima-Crystal}
J.~Beck and H.~Nakajima.
\newblock Crystal bases and two-sided cells of quantum affine algebras.
\newblock {\em Duke Math. J.}, 123(2):335--402, 2002.

\bibitem{Chen-Rootvector}
X.~Chen.
\newblock Root vectors of the composition algebra of the {K}ronecker algebra.
\newblock {\em Algebra Discrete Math.}, 2004(1):37--56, 2004.

\bibitem{DengDuXiao2007Generic}
B.~Deng, J.~Du, and J.~Xiao.
\newblock Generic extensions and canonical bases for cyclic quivers.
\newblock {\em Canad. J. Math.}, 59(6):1260--1283, 2007.

\bibitem{green1995hall}
J.~A. Green.
\newblock Hall algebras, hereditary algebras and quantum groups.
\newblock {\em Invent. Math.}, 120(1):361--377, 1995.

\bibitem{Li-yiqiang_canonical_monomial}
Y.~Li.
\newblock Notes on affine canonical and monomial bases.
\newblock {\em arXiv preprint arXiv:math/0610449}, 2006.

\bibitem{Lin_Xiao_Zhang_Representations_of_tame_quivers_and_affine_canonical_bases}
Z.~Lin, J.~Xiao, and G.~Zhang.
\newblock Representations of tame quivers and affine canonical bases.
\newblock {\em Publ. Res. Inst. Math. Sci.}, 47:825--885, 2011.

\bibitem{Lusztig_Canonical_bases_arising_from_quantized_enveloping_algebra}
G.~Lusztig.
\newblock Canonical bases arising from quantized enveloping algebras.
\newblock {\em J. Amer. Math. Soc.}, 3(2):447--498, 1990.

\bibitem{Lusztig_Quivers_perverse_sheaves_and_the_quantized_enveloping_algebras}
G.~Lusztig.
\newblock Quivers, perverse sheaves, and quantized enveloping algebras.
\newblock {\em J. Amer. Math. Soc.}, 4(2):365--421, 1991.

\bibitem{Lusztig_Affine_quivers_and_canonical_bases}
G.~Lusztig.
\newblock Affine quivers and canonical bases.
\newblock {\em Publications Math{\'e}matiques de l'IH{\'E}S}, 76:111--163,
  1992.

\bibitem{Lusztig_Introduction_to_quantum_groups}
G.~Lusztig.
\newblock {\em Introduction to quantum groups}.
\newblock Springer, 2010.

\bibitem{Mcgerty2004The_Kronecker_quiver}
K.~Mcgerty.
\newblock The {K}ronecker quiver and bases of quantum affine $\mathfrak{sl}_2$.
\newblock {\em Adv. Math.}, 197(2):411--429, 2005.

\bibitem{Nakajima_Crystal_canonical_and_PBW_bases_of_quantum_affine_algebras}
H.~Nakajima.
\newblock Crystal, canonical and {P}{B}{W} bases of quantum affine algebras.
\newblock In {\em Algebraic groups and homogeneous spaces}, pages 389--421.
  Narosa Publishing House, New Delhi, India, 2007.

\bibitem{Ringel_Hall_algebras_and_quantum_groups}
C.~M. Ringel.
\newblock Hall algebras and quantum groups.
\newblock {\em Invent. Math.}, 101(1):583--591, 1990.

\bibitem{Rongel_The_Hall_algebra_approach_to_quantum_groups}
C.~M. Ringel.
\newblock The {H}all algebra approach to quantum groups.
\newblock In {\em Proceedings E.L.A.M.}, pages 85--114. Aportaciones
  Matem{\'a}ticas Comunicaciones, 1995.

\bibitem{Xiao-Xu-Zhao}
J.~Xiao, H.~Xu, and M.~Zhao.
\newblock Tame quivers and affine bases {I}: a {H}all algebra approach to the
  canonical bases.
\newblock {\em arXiv preprint arXiv:2107.03095}, 2021.

\bibitem{Zhang2000PBW}
P.~Zhang.
\newblock {PBW}-basis for the composition algebra of the {K}ronecker algebra.
\newblock {\em J. Reine Angew. Math.}, 2000(527):97--116, 2000.

\end{thebibliography}

\end{document}